\documentclass[reqno, 11pt]{amsart}

\usepackage{amsfonts, amsthm, amsmath, amssymb,bbold,dsfont}
\usepackage{hyperref}
\usepackage{tikz}
\hypersetup{colorlinks=false}

\usepackage{bbm}  

 \usepackage{tabu}

\usepackage[margin=1in]{geometry}

\usepackage{helvet}

\RequirePackage{mathrsfs} \let\mathcal\mathscr

\usepackage{hyperref}
\usepackage{dsfont}
\usepackage{upgreek}
\usepackage{mathabx,yfonts}
\usepackage{enumerate}

\numberwithin{equation}{section}

\newtheorem{theorem}{Theorem}[section] 
\newtheorem{lemma}[theorem]{Lemma}

\theoremstyle{definition}

 \newtheorem*{acknowledgements}{Acknowledgements}
\newtheorem{remark}[theorem]{Remark}

 \newtheorem*{notation}{Notation}

\renewcommand{\phi}{\varphi}

\renewcommand{\leq}{\leqslant}

\renewcommand{\geq}{\geqslant}

\newcommand{\M}{\mathbf{M}}

\newcommand{\z}{\mathbf{z}}

\DeclareMathOperator{\nm}{Nm}
\DeclareMathOperator{\n}{N}

\DeclareSymbolFont{bbold}{U}{bbold}{m}{n}
\DeclareSymbolFontAlphabet{\mathbbold}{bbold}

\renewcommand{\epsilon}{\varepsilon}

\renewcommand{\leq}{\leqslant}
\renewcommand{\geq}{\geqslant}

\title[General bilinear forms in the Jacobi symbol over hyperbolic regions]{General bilinear forms in the Jacobi symbol over hyperbolic regions}

\author{Cameron Wilson}
\address{Department of Mathematics\\
University of Glasgow  \\ G12~8QQ United Kingdom}
\email{c.wilson.6@research.gla.ac.uk}

\date{\today}

\begin{document}
\begin{abstract} We study averages involving the Jacobi quadratic symbol $(\frac{n}{m})$ in regions where the product $mn$ is bounded by a large parameter. We show that these averages exhibit cancellation whenever the summation is restricted to square-free integers bounded away from the axes.
\end{abstract}

\maketitle
\setcounter{tocdepth}{1}
\tableofcontents
\section{Introduction}\label{s:intro}
General bilinear forms in the Jacobi symbol over rectangular regions have been a focus of intensive study in recent decades (see Heath-Brown \cite{Heath-Brown} and Friedlander--Iwaniec \cite{F--I}). They have been used in many problems, for example:
\begin{itemize}
\item Values of $L$-functions (Soundararajan  \cite{S});
\item 4-ranks of class groups (Fouvry--Kl\"{u}ners \cite{F--K});
\item Manin's conjecture (Browning--Heath-Brown  \cite{TDBHB});
\item Bateman--Horn's conjecture on average (Baier--Zhao \cite{BSZL}).
\end{itemize}
Problems of current interest in arithmetic geometry, however, require a study of such bilinear forms with height conditions of a geometric nature. In particular, versions where the sums over rectangular regions are replaced by those over hyperbolic regions appear when considering the local solubility of quadrics over surfaces parameterised by $\mathbb{P}^1\times\mathbb{P}^1$ (see the top paragraph of
\cite[pp. 3]{BLS}). These sums are of the form
\begin{equation}    \label{typicalhyperbolicsum}
    \sum_{\substack{n,m\in\mathbb{N}\\1\leq nm \leq T}} a_n b_m \Big(\frac{n}{m}\Big),
\end{equation}
where $(\frac{n}{m})$ is the Jacobi quadratic symbol and $(a_n)$ and $(b_m)$ are arbitrary complex sequences with $\lvert a_n\rvert,\lvert b_m\rvert\leq 1$. For general choices of complex sequences $a_n$ and $b_m$ these sums do not give much cancellation - for example, we will see later that
\begin{equation*}
    \sum_{\substack{n,m\in\mathbb{N}\\2\nmid nm\\ 1\leq nm\leq T}} \Big(\frac{n}{m}\Big)\gg T,
\end{equation*}
which gives only logarithmic saving over the hyperbolic region of area $T(\log T)$. The main contribution of this sum will be seen to come from the points where either $n$ or $m$ is a square. This contribution is explained by the fact that such points have relatively large density in the hyperbolic region compared to their density in a rectangular one. For this reason we turn to study sums over pairs $(n,m)$ where $n$ and $m$ are odd and square-free.\\

\notation Throughout this paper $\sideset{}{^*}\sum$ will denote a sum over odd, square-free integers. As usual, $\mu$ will denote the M\"obius function.\\

In this case, however, there may still be a large contribution from points close to the axes, particularly from the lines $n=1$ and $m=1$. Another example which gives very small cancellation is the following: choosing $a_n = (\frac{n}{11})$ and $b_m$ to be the characteristic function for the condition $m=11$, one sees that
$$\sideset{}{^*}\sum_{\substack{1<n,m\leq T\\1\leq nm\leq T}}\mu^2(2nm)a_nb_m\Big(\frac{n}{m}\Big) = \sum_{1<n\leq T/11}\mu^2(22n)\gg T.$$
It is therefore clear that in order to obtain cancellation, we must impose the extra condition that $n,m>z$ for some parameter $z=z(T)$ which tends to infinity with $T$. We will show that such restrictions will give appropriate cancellation.
\begin{theorem}\label{first}
Let $T,z\geq 2$ and let $(a_n)$, $(b_m)$ be any complex sequences such that $\lvert a_n\rvert,\lvert b_m\rvert \leq 1$. If there exists an $\epsilon>0$ such that $z\geq T^{1/3-\epsilon}$, then
\[
\sideset{}{^*}\sum_{\substack{z < n,m \leq T\\ nm \leq T}} a_n b_m \Big(\frac{n}{m}\Big) \ll_{\epsilon} \frac{T^{1+\epsilon}}{z^{1/2}},
\]
where the implied constant depends at most on $\epsilon$. If there exists an $\epsilon>0$ such that $z\leq T^{1/3-\epsilon}$, then
\[
\sideset{}{^*}\sum_{\substack{z < n,m \leq T\\ nm \leq T}} a_n b_m \Big(\frac{n}{m}\Big) \ll_{\epsilon} \frac{T(\log T)^3}{z^{1/2}},
\]
where the implied constant depends at most on $\epsilon$.
\end{theorem}
Theorem \ref{first} fails to give a saving over the trivial bound of $T(\log T)$ if $z \ll (\log T)^{4}$. This is satisfactory for most applications, however it is also possible to obtain cancellation for any $z$ which tends to infinity with $T$ with the cost of a smaller exponent of $z$.
\begin{theorem}\label{main}
For all $T,z\geq 2$ and all complex sequences $(a_n)$, $(b_m)$ such that $\lvert a_n\rvert,\lvert b_m\rvert \leq 1$ we have,
\[
    \sideset{}{^*}\sum_{\substack{z < n,m \leq T\\ nm \leq T}} a_n b_m \Big(\frac{n}{m}\Big) \ll \frac{T(\log T)}{z^{1/4}},
\]
where the implied constant is absolute.\\
\end{theorem}
\remark These results show that the majority of the sum over the hyperbolic region will contribute to an error term. Indeed, if $z\leq T^{1/3-\epsilon}$ for some $\epsilon>0$, we may write
\begin{align*}
\sideset{}{^*}\sum_{\substack{1\leq n,m\leq T\\1\leq nm \leq T}} a_n b_m \Big(\frac{n}{m}\Big) &= \sideset{}{^*}\sum_{\substack{1\leq n\leq z\\1\leq m \leq T/n}} a_n b_m \Big(\frac{n}{m}\Big) + \sideset{}{^*}\sum_{\substack{1\leq m\leq z\\ 1\leq n \leq T/m}} a_n b_m \Big(\frac{n}{m}\Big)\\
&+ \sideset{}{^*}\sum_{\substack{z< n,m\leq T\\1\leq nm \leq T}} a_n b_m \Big(\frac{n}{m}\Big) - \sideset{}{^*}\sum_{\substack{1\leq n\leq z\\1\leq m\leq z}}a_nb_m\Big(\frac{n}{m}\Big).
\end{align*}
Then the third sum on the right hand side may be bounded by $\frac{T(\log T)^3}{z^{1/2}}$ using Theorem \ref{first} and the fourth sum may be bounded by $z^{3/2}$ using a well known result of Heath-Brown \cite{Heath-Brown} (see equation \eqref{HB} below). These error terms will be sufficiently small whenever $z = (\log T)^A$ for $A>6$. Choosing $z$ to be as such, we may deal with the remaining two sums using familiar methods for averages of characters with small modulus such as the tools used in the Siegel--Walfisz theorem. For example, see \cite[equation (24)]{F--I}.\\

\remark It is impossible to improve the exponent of $z$ in Theorem \ref{main} to be $>1$ in general: using a similar example to before, if we take any $z\leq T^{1/4}$, $p$ a prime satisfying $z<p\leq 2z$, $a_n = (\frac{n}{p})$ and $b_m$ the characteristic function for the condition $m=p$, one finds that
$$\sideset{}{^*}\sum_{\substack{z<n,m\leq T\\1\leq nm\leq T}}a_nb_m\Big(\frac{n}{m}\Big) = \sum_{z<n\leq T/p}\mu^2(2pn)= \frac{2}{3(1+1/p)\zeta(2)}\left(\frac{T}{p}-z\right)+O\left(\frac{T^{1/2}}{p^{1/2}}\right).$$
This is $\gg \frac{T}{p}\gg \frac{T}{z}$ since $2pz < 4z^2\leq T$. If this were $O(\frac{T(\log T)}{z^{\alpha}})$, then $z^{\alpha-1} = O(\log T)$, which may be contradicted by taking $z=(\log T)^A$ for $A>0$ suitably large.\\

In order to prove Theorems \ref{first} and \ref{main} we will actually prove more general results. In particular we will prove the following:
\begin{theorem}\label{general1}
Let $T,z\geq 2$, $c\geq 0$, and let $(a_n)$, $(b_m)$ be any complex sequences such that $\lvert a_n\rvert,\lvert b_m\rvert \leq 1$. If there exists an $\epsilon>0$ such that $z\geq T^{1/3-\epsilon}$ then,
    \[
    \sideset{}{^*}\sum_{\substack{z < n,m \leq T\\ nm \leq T}} (nm)^c a_n b_m \Big(\frac{n}{m}\Big) \ll_{c,\epsilon} \frac{T^{1+c+\epsilon}}{z^{1/2}},
    \]
where the implied constant depends at most on $c$ and $\epsilon$. If there exists an $\epsilon>0$ such that $z\leq T^{1/3-\epsilon}$, then
\[
\sideset{}{^*}\sum_{\substack{z < n,m \leq T\\ nm \leq T}} (nm)^c a_n b_m \Big(\frac{n}{m}\Big) \ll_{c,\epsilon} \frac{T^{1+c}(\log T)^3}{z^{1/2}},
\]
where the implied constant depends at most on $c$ and $\epsilon$.
\end{theorem}

\begin{theorem}\label{general2}
For all $T,z\geq 2$, $c\geq 0$, and all complex sequences $(a_n)$, $(b_m)$ such that $\lvert a_n\rvert,\lvert b_m\rvert \leq 1$. Then
    \[
    \sideset{}{^*}\sum_{\substack{z < n,m \leq T\\ nm \leq T}} (nm)^c a_n b_m \Big(\frac{n}{m}\Big) \ll_c \frac{T^{1+c}(\log T)}{z^{1/4}},
    \]
where the implied constant depends at most on $c$.
\end{theorem}
Theorems \ref{first} and \ref{main} will then follow from the cases where $c=0$. The methods in \cite{F--I,Heath-Brown} do not interact well with the hyperbolic region as they exploit the linear structure of the rectangular regions through the use of H\"older's inequality. We will circumvent this problem by applying the following version of Perron's formula to eliminate the hyperbolic height condition \cite[Lemma 2.2]{Harman}:
\begin{equation}\label{Perron}
\frac{1}{\pi}\int_{-R}^{R} (nm)^{it} f_{\tau}(t) dt = \mathds{1}(nm\leq \tau) + O(R^{-1}\lvert\log(nm)-\log(\tau)\rvert^{-1}),
\end{equation}
where $f_{\tau}(t)=\frac{\sin(t\log(\tau))}{t}$ and 
\begin{equation*}
    \mathds{1}(\mu\leq\tau) = \begin{cases}
        1\;\text{if}\;\mu\leq\tau
        \\ 0\;\text{if}\;\mu>\tau.
    \end{cases}
\end{equation*} 
This will allow us to apply existing results. In particular, we apply Corollary $4$ and Theorem $1$ of \cite{Heath-Brown}:
\begin{equation}\label{HB}
    \sum_{\substack{m\leq M\\2\nmid m}}\sum_{n\leq N}a_nb_m\Big(\frac{n}{m}\Big) \ll_{\epsilon} (MN)^{\epsilon}(MN^{1/2}+M^{1/2}N),
\end{equation}
and, for $I\subseteq [1,N]\cap\mathbb{N}$ of size $\lvert I\rvert$,
\begin{equation} \label{HB1}
        \sideset{}{^*}\sum_{m\leq M}\Big\lvert\sideset{}{^*}\sum_{n\in I} a_n\Big(\frac{n}{m}\Big)\Big\rvert^2\ll_{\epsilon} (MN)^{\epsilon}(\max(M,N))\lvert I\rvert,
\end{equation}
for any $\epsilon >0$. It will be also be necessary to make use of the following result:
\begin{equation} \label{Elliot}
        \sideset{}{^*}\sum_{m\leq M}\Big\lvert\sideset{}{^*}\sum_{n\in I} a_n\Big(\frac{n}{m}\Big)\Big\rvert^2\ll (M+N^2\log(N))\lvert I\rvert.
\end{equation}
This inequality goes back to Elliot \cite{Elliot}, but was proven by Heath-Brown \cite[Equation (6)]{Heath-Brown}. If $N \leq M^{1/2}$ then this bound is superior to \eqref{HB1}. In our proof, this will be necessary when lopsided rectangles appear in our coverings of the hyperbolic region. We cannot use \eqref{HB1} for such lopsided rectangles, since if $z = (\log T)^A$ for some $A>0$, we will obtain bounds of the form $\ll\frac{T^{1+c+\epsilon}}{(\log T)^{A/2}}$, which just fails to give Theorem \ref{general1}.\\

To prove Theorem \ref{general1} we will use these results along with a diadic covering of the hyperbolic regions. In order to obtain saving for arbitrary $z$ in Theorem \ref{general2} it will be necessary to cover parts of the hyperbolic region with rectangles of equal width before applying Cauchy--Schwarz and \eqref{Elliot} to each of these rectangles and summing over the results.\\

Lastly, we prove asymptotics for the Jacobi sums over all odd integers within the hyperbolic region, namely we have the following:
\begin{theorem}\label{asymptotic}
For all $T\geq 2$,
\[
\sum_{\substack{1\leq n,m\leq T\\2\nmid nm\\nm\leq T}} \Big(\frac{n}{m}\Big) = \left(\frac{6\zeta(2)}{7\zeta(3)}\right)T + O(T^{3/4}(\log T))
\]
where $\zeta$ is the Riemann-zeta function.
\end{theorem}
This result is obtained using Dirichlet's hyperbola method. A similar method will give
\[
\sideset{}{^*}\sum_{\substack{1\leq n,m\leq T\\nm\leq T}} \Big(\frac{n}{m}\Big) \sim c'T.
\]
for some constant $c'>0$, proving that we may not obtain saving in the square-free setting when we include points close to the axes.

\acknowledgements We are very grateful to the referee for reading the paper carefully and providing many helpful comments, including the suggestion to use Perron's formula which greatly streamlined parts of the argument throughout the paper. The author would also like to thank his supervisor, Efthymios Sofos, for helping to shape this paper to its current form, as well as the Carnegie Trust for the Universities of Scotland for the funding received through their PhD Scholarship Programme.

\section{Proof of Theorem \ref{general1}}\label{gen1proof}
Our strategy for this proof will be to cover the hyperbolic region in diadic rectangles and apply \eqref{Perron}-\eqref{Elliot}. Set
\[
S(T) = \sideset{}{^*}\sum_{\substack{z < n,m \leq T\\ nm \leq T}} (nm)^c a_n b_m \Big(\frac{n}{m}\Big).
\]
If $z>T^{1/2}$ then the sum is $0$, so that the bound is trivially true. If $z = T^{1/2}$ then the sum has magnitude $\leq 1$, so again the bound is trivial. We are now left with $z<T^{1/2}$. Our first step is to split the sum over the hyperbolic region into $4$ pieces. We write
\[
S(T) = S_1(T) + S_2(T) + S_3(T) - S_4(T)
\]
where
\begin{align*}
    S_1(T)&=\sideset{}{^*}\sum_{\substack{z_1 < n,m \leq T\\nm\leq T}} (nm)^ca_n b_m \Big(\frac{n}{m}\Big);\\
    S_2(T)&=\sideset{}{^*}\sum_{z < n \leq z_1\;}\sideset{}{^*}\sum_{\;z < m\leq \frac{T}{n}} (nm)^c a_n b_m \Big(\frac{n}{m}\Big);\\
    S_3(T)&= \sideset{}{^*}\sum_{z < m \leq z_1\;}\sideset{}{^*}\sum_{\;z < n\leq \frac{T}{m}} (nm)^c a_n b_m \Big(\frac{n}{m}\Big);\\
    S_4(T)&=\sideset{}{^*}\sum_{z < n \leq z_1\;}\sideset{}{^*}\sum_{\;z < m\leq z_1} (nm)^c a_n b_m \Big(\frac{n}{m}\Big).
\end{align*}
where $z_1 = \max(z,\frac{T^{1/3}}{\log T})$. We now aim to bound each of these sums individually. First note that if $z> \frac{T^{1/3}}{\log T}$, $S_2(T)=S_3(T)=S_4(T)=0$ and so we only need to consider them whenever $z\leq\frac{T^{1/3}}{\log T}$. For $S_4(T)$ we apply \eqref{HB} with $N=M=\frac{T^{1/3}}{\log T}$: divide and multiply the sum by $T^c$ so that we have $\lvert\frac{n^c}{T^{c/2}}a_n\rvert$, $\lvert\frac{m^c}{T^{c/2}}b_m\rvert \leq 1$. Then \eqref{HB} with $\epsilon < 1/6$ gives
\begin{align*}
S_4(T) \ll_c \begin{cases} T^{2/3+c}&\;\text{if}\;z\leq \frac{T^{1/3}}{\log T},\\ 0 &\;\text{if}\;z> \frac{T^{1/3}}{\log T}.\end{cases}
\end{align*}
We now turn to the remaining $3$ sums. If $z\leq \frac{T^{1/3}}{\log T}$, then $S_2(T)$ and $S_3(T)$ may be dealt with using symmetric arguments as a consequence of reciprocity for Jacobi symbols, $(\frac{n}{m}) = (-1)^{\frac{(n-1)(m-1)}{4}}(\frac{m}{n})$, and so we only need to deal with one: since $n$ and $m$ are odd and square-free, we may split $S_3(T)$ into $4$ sums using the conditions $n,m\equiv1\;\text{or}\;3\pmod{4}$ and then reciprocity will give $4$ sums in the same form as $S_2(T)$. Thus we only need to consider $S_1(T)$ and $S_2(T)$. We aim to use Perron's formula. For $S_3(T)$ we split $(z_1,T]$ into diadic intervals to obtain $\ll(\log T)^2$ diadic regions $(N,2N]\times(M,2M]$ where $N,M\in(z_1,T]$. For $S_2(T)$, split the intervals $(z,z_1]$ and $(z,T]$ into diadic intervals to obtain $\ll(\log T)^2$ diadic regions $(N,2N]\times(M,2M]$ where $N\in(z,z_1]$ and $M\in(z,T]$. This will give bounds of the form
\begin{align*}
S_1(T) &\ll (\log T)^2\max_{\substack{z_1<N\leq T\\z_1<M\leq T\\ NM\leq T}}\lvert S(T;N,M)\rvert;\\
S_2(T) &\ll (\log T)^2\max_{\substack{z<N\leq z_1\\z<M\leq T\\ NM\leq T}}\lvert S(T;N,M)\rvert;
\end{align*}
where in each case
\[
S(T;N,M) = \sideset{}{^*}\sum_{\substack{N < n \leq 2N\\M < m \leq 2M \\ nm \leq T}} (nm)^c a_n b_m \Big(\frac{n}{m}\Big).
\]
Next, we apply Perron's formula to deal with the hyperbolic conditions. 
Let $\theta\in[-1/2,1/2]$ be such that $T+\theta\in\mathbb{Z}+\frac{1}{2}$ and take $\tau = T+\theta$ in \eqref{Perron}. Then \eqref{Perron} becomes
\[
\frac{1}{\pi}\int_{-R}^{R} (nm)^{it} f_{T+\theta}(t) dt = \mathds{1}(nm\leq T) + O(R^{-1}\lvert\log(nm)-\log(T+\theta)\rvert^{-1})
\]
for any $R>0$ where $f_{T+\theta}(t)=\frac{\sin(t\log(T+\theta))}{t}$. Noting that here $(\log(nm)-\log(T+\theta))\gg\frac{1}{T}$, we substitute this into $S(T;N,M)$ to obtain
\[
S(T;N,M) = \frac{1}{\pi} \int_{-R}^{R} f_{T+\theta}(t)\sideset{}{^*}\sum_{\substack{N < n \leq 2N\\M < m \leq 2M }} (nm)^{c+it} a_n b_m \Big(\frac{n}{m}\Big) dt + O\left(\frac{(NM)^{1+c}T}{R}\right).
\]
Before we move forward, we deal with the $(nm)^{c+it}$ term: write $1 = \frac{(4NM)^c}{(4NM)^c}$, then by setting $\widetilde{a}_n = \frac{n^{c}}{(2N)^c}a_n$ and $\widetilde{b}_m = \frac{m^{c}}{(2M)^c}b_m$ we have
\[
\sideset{}{^*}\sum_{\substack{N < n \leq 2N\\M < m \leq 2M }} (nm)^{c+it} a_n b_m \Big(\frac{n}{m}\Big) = (4NM)^c \sideset{}{^*}\sum_{\substack{N < n \leq 2N\\M < m \leq 2M }}  n^{it}\widetilde{a}_n m^{it}\widetilde{b}_m \Big(\frac{n}{m}\Big)
\]
where $\lvert n^{it}\widetilde{a}_n\rvert$, $\lvert m^{it}\widetilde{b}_m\rvert \leq 1$. By applying \eqref{HB} and substituting into $S(T;N,M)$ we get
\begin{equation}\label{HBYield}
    S(T;N,M) \ll_{\epsilon} 4^c (MN)^c \int_{-R}^{R} \lvert f_{T+\theta}(t)\rvert dt (MN)^{\epsilon}\left(MN^{1/2}+M^{1/2}N\right) + \frac{(NM)^{1+c}T}{R}.
\end{equation}
By instead applying Cauchy--Schwarz and \eqref{Elliot} we obtain
\begin{equation}\label{ElliotYield}
    S(T;N,M) \ll_{\epsilon} 4^c (MN)^c \int_{-R}^{R} \lvert f_{T+\theta}(t)\rvert dt \left(MN^{1/2}+M^{1/2}N^{3/2}(\log N)^{1/2}\right) + \frac{(NM)^{1+c}T}{R}.
\end{equation}
We will apply \eqref{HBYield} to the diadic regions in $S_1(T)$ to obtain
\begin{align*}
S_1(T) &\ll_{c,\epsilon} T^{c}(\log T)^2 \int_{-R}^{R}\lvert f_{T+\theta}(t)\rvert dt \max_{\substack{z_1<N\leq T\\z_1<M\leq T\\ NM\leq T}}(MN)^{\epsilon}\left(MN^{1/2}+M^{1/2}N\right) + \frac{T^{2+c}(\log T)^2}{R}\\
&\ll_{c,\epsilon} T^{c}(\log T)^2 \int_{-R}^{R}\lvert f_{T+\theta}(t)\rvert dt \left(\frac{T^{1+\epsilon}}{z_1^{1/2}}\right) + \frac{T^{2+c}(\log T)^2}{R}.
\end{align*}
For $S_2(T)$ however, we may assume $z\leq\frac{T^{1/3}}{\log T}$ so that the diadic rectangles either have the lopsided condition $N(\log N)^{1/2}\leq M^{1/2}\leq \frac{T^{1/2}}{z^{1/2}}$ or have $M^{1/2}\leq N(\log N)^{1/2}\ll \frac{T^{1/3}}{(\log T)^{1/2}}$. Thus we may use \eqref{ElliotYield} to obtain
\begin{align*}
S_2(T)&\ll_{c} T^c(\log T)^2 \int_{-R}^{R} \lvert f_{T+\theta}(t)\rvert dt \max_{\substack{z<N\leq \frac{T^{1/3}}{\log T}\\z<M\leq T\\ NM\leq T}}(MN^{1/2}+M^{1/2}N^{3/2}(\log N)^{1/2}) + \frac{T^{2+c}(\log T)^2}{R}\\
&\ll_{c} T^c(\log T)^2 \int_{-R}^{R} \lvert f_{T+\theta}(t)\rvert dt \left(T^{1/2}\left(\frac{T^{1/2}}{z^{1/2}}+\frac{T^{1/3}}{(\log T)^{1/2}}\right)\right) + \frac{T^{2+c}(\log T)^2}{R}\\
&\ll_{c} T^c(\log T)^2 \int_{-R}^{R} \lvert f_{T+\theta}(t)\rvert dt\left(\frac{T}{z^{1/2}}\right) + \frac{T^{2+c}(\log T)^2}{R},
\end{align*}
if $z\leq \frac{T^{1/3}}{\log T}$ and $S_2(T)=0$ otherwise. Choosing $R = T^2(\log T)^2$ the integral becomes bounded by $(\log T)$ therefore giving
\[
S_1(T)\ll_{c,\epsilon}\frac{T^{1+c+\epsilon}(\log T)^3}{z_1^{1/2}},
\]
and
\begin{align*}
S_2(T)\ll_{c} \begin{cases} \frac{T^{1+c}(\log T)^3}{z^{1/2}}&\text{if}\;z\leq \frac{T^{1/3}}{\log T},\\
\;\;\;\;\;\;0&\text{if}\;z>\frac{T^{1/3}}{\log T}.
\end{cases}
\end{align*}
Finally, recalling that $S_2(T)$ and $S_3(T)$ are symmetrically equivalent and that $z_1=\max(z,\frac{T^{1/3}}{\log T})$, we may put all of these bounds together to obtain
\begin{align*}
S(T)\ll_{c,\epsilon} \begin{cases} \frac{T^{1+c}(\log T)^3}{z^{1/2}}+T^{5/6+c+\epsilon}(\log T)^{7/2}+T^{2/3+c}&\text{if}\;z\leq \frac{T^{1/3}}{\log T},\\
\frac{T^{1+c+\epsilon}(\log T)^3}{z^{1/2}}&\text{if}\;z>\frac{T^{1/3}}{\log T},
\end{cases}
\end{align*}
for any $\epsilon>0$. Now suppose there exists an $\epsilon>0$ such that $z\geq T^{1/3-\epsilon}$. Then if $z>\frac{T^{1/3}}{\log T}$ we use the second case of the above bound with $\epsilon/2$ to obtain:
\[
S(T)\ll_{c,\epsilon} \frac{T^{1+c+\epsilon/2}(\log T)^3}{z^{1/2}}\ll_{c,\epsilon} \frac{T^{1+c+\epsilon}}{z^{1/2}}.
\]
If $T^{1/3-\epsilon}\leq z\leq \frac{T^{1/3}}{\log T}$ then consider the first bound with $\epsilon/2$. Then,
\[
S(T)\ll_{c,\epsilon} \frac{T^{1+c}(\log T)^3}{z^{1/2}}+T^{5/6+c+\epsilon/2}(\log T)^{7/2}+T^{2/3+c} \ll_{c,\epsilon} \frac{T^{1+c+\epsilon}}{z^{1/2}}.
\]
Lastly, if there exists an $\epsilon>0$ such that $z\leq T^{1/3-\epsilon}$, then consider the bound for $z\leq \frac{T^{1/3}}{\log T}$ with $\epsilon/2$. In this case the first term dominates since
\[
\frac{T^{1+c}(\log T)^3}{z^{1/2}}\geq T^{5/6+c+\epsilon}(\log T)^3 \gg T^{5/6+c+\epsilon/2}(\log T)^{7/2},
\]
which implies the result.

\section{Proof of Theorem \ref{general2}}\label{gen2proof}
The key idea of this proof is to cover parts of the hyperbolic region with rectangles of equal width and apply Theorem \ref{general1} along with Cauchy--Schwarz and \eqref{Elliot} to the sums over each of these rectangle and then sum over the results. The following lemma encodes the covering we will use:

\begin{lemma}\label{cover}
Fix $c\geq 0$ and $0<\delta\leq 1/2$. Then for any $T\geq 2$ and any $2\leq z < T^{\delta}$ we have,
\begin{align*}
\sideset{}{^*}\sum_{z < n \leq T^{\delta}\;}\sideset{}{^*}\sum_{\;z < m\leq \frac{T}{n}} (nm)^c a_n b_m \Big(\frac{n}{m}\Big) = \sum_{k\in H}\;\;\sideset{}{^*}\sum_{n\in I_k}\;\;\sideset{}{^*}\sum_{m\in J_k} (nm)^c a_n b_m \Big(\frac{n}{m}\Big) + O\left(\frac{T^{1+c}}{z^{1/2}}\right),
\end{align*}
where
\begin{align*}
&H=\left\{k\in\mathbb{N}:z^{1/2}\leq k\leq \frac{T^{\delta}}{z^{1/2}}-1\right\};\\
&I_k = \left\{n\in\mathbb{N}: z^{1/2}k<n\leq z^{1/2}(k+1)\right\};\\
&J_k = \left\{m\in\mathbb{N}: z<m\leq \frac{T}{z^{1/2}(k+1)}\right\}.
\end{align*}
\end{lemma}

\remark By partitioning the interval over $n$ into intervals of equal length we are then able to make use of the fact that, close to the hyperbolic curve, the gradient $\frac{\mathrm{d}(x/n)}{\mathrm{d}n}=-\frac{x}{n^2}$, decreases rapidly in magnitude. This observation will lead to the regions leftover from the covering boxes having small volume.

\begin{proof}
We begin by partitioning the interval $(0,T^{\delta}]$ into $\frac{T^{\delta}}{z^{1/2}}$ intervals of equal length $z^{1/2}$, say $(kz^{1/2},(k+1)z^{1/2}]$ for integers $0\leq k\leq \frac{T^{\delta}}{z^{1/2}}-1$, and intersecting this partition with $(z,T^{\delta}]$. Then notice that
$$\bigcup_{k\in H} I_k \subseteq (z,T^{\delta}]\cap\mathbb{N},$$
where the leftover part of this partition, $L'$, satisfies
$$L' = ((z,T^{\delta}]\cap\mathbb{N})\setminus \bigcup_{k\in H} I_k \subseteq ((z,z+z^{1/2}]\cup(T^{\delta}-z^{1/2}, T^{\delta}])\cap\mathbb{N}.$$
Fix a $k$. Then for an $n\in I_k$, the summation index $m$ ranges from $z$ to $\frac{T}{n}$ where $\frac{T}{(k+1)z^{1/2}}\leq \frac{T}{n} < \frac{T}{kz^{1/2}}$. To create our rectangles we split all ranges over $m$ into a range $z<m\leq\frac{T}{(k+1)z^{1/2}}$, giving us the intervals $J_k$, and $\frac{T}{(k+1)z^{1/2}}<m\leq \frac{T}{n}$. Notice that for $n\in I_k$, $(\frac{T}{(k+1)z^{1/2}},\frac{T}{n}]\cap\mathbb{N} \subseteq J'_k = \{m\in\mathbb{N}: \frac{T}{z^{1/2}(k+1)}<m\leq \frac{T}{z^{1/2}k}\}$. Combining the ranges, for each $k$ we have a rectangle $I_k\times J_k$ and a small section $L_k = \{(n,m)\in\mathbb{N}:n\in I_k,\; \frac{T}{(k+1)z^{1/2}}<m\leq \frac{T}{n}\}$ close to the hyperbolic curve which is contained in the small rectangle $I_k\times J'_k$. We also have the leftover regions coming from $n\in L'$, $L = \{(n,m)\in\mathbb{N}^2:n\in L',\; z<m\leq\frac{T}{n}\}$. Then we have:
\begin{align*}
\{(n,m)\in\mathbb{N}^2:nm\leq T,\; n\leq T^{\delta},\;n,m>z\} = \bigcup_{k\in H}(I_k\times J_k) \cup \bigcup_{k\in H} L_k \cup L.
\end{align*}
See Figure $1$ for an illustration of these sets.
\begin{figure}\label{coverfigure}
\begin{tikzpicture}
\draw[->,thick] (0,-1)--(0,10);
\draw[->,thick] (-1,0)--(13.33,0);
\draw[thick,scale=1,domain=4:13.33,smooth,variable=\t] plot ({\t},{40/\t});
\draw[dotted,thick] (0,2.5)--(2.5,2.5);
\draw[dashed] (2.5,2.5)--(13.33,2.5);
\draw[dotted,thick] (0,7)--(13.33,7);
\draw[dashed] (2.5,4)--(10,4);
\draw[dashed] (2.5,5)--(8,5);
\draw[dashed] (2.5,4)--(2.5,5);
\draw[dotted,thick] (0,4)--(2.5,4);
\draw[dotted,thick] (0,5)--(2.5,5);
\draw[dotted,thick] (8,5)--(10,5);
\draw[dotted,thick] (10,4)--(10,5);
\draw[dashed] (8,4)--(8,5);
\draw[dotted,thick] (8,0)--(8,4);
\draw[dotted,thick] (10,0)--(10,4);
\draw[dotted,thick] (2.5,0)--(2.5,10);
\draw[dashed] (2.5,2.5)--(2.5,3);
\draw[dotted,thick] (0,3)--(2.5,3);
\draw[dashed] (2.5,3)--(13.33,3);

\node[] at (2.5,-1) {$z$};
\node[] at (-1,7) {$T^{\delta}$};
\node[] at (-1,6.2) {$\vdots$};
\node[] at (-1,2.5) {$z$};
\node[] at (-1,3.3) {$\vdots$};
\node[] at (-1,4) {$kz^{1/2}$};
\node[] at (-1,5) {$(k+1)z^{1/2}$};
\node[] at (5,-1) {$\cdots$};
\node[] at (8,-1) {$\frac{T}{(k+1)z^{1/2}}$};
\node[] at (10,-1) {$\frac{T}{kz^{1/2}}$};
\node[] at (-1,10) {$m$};
\node[] at (15,-1) {$n$};
\node[] at (5,4.5) {$I_k\times J_k$};
\node[] at (8.5,4.25) {$L_k$};
\node[] at (8.25,2.75) {$L$};
\end{tikzpicture}
\caption{Lemma \ref{cover} Illustration}
\end{figure}

It follows that
\begin{align*}
    \sideset{}{^*}\sum_{z < n \leq T^{\delta}\;}\sideset{}{^*}\sum_{\;z < m\leq \frac{T}{n}} (nm)^c a_n b_m \Big(\frac{n}{m}\Big) = &\sum_{k\in H} \;\sideset{}{^*}\sum_{n\in I_k}\;\sideset{}{^*}\sum_{m\in J_k}(nm)^c a_nb_m\Big(\frac{n}{m}\Big)\\ +&\sum_{k\in H}\;\sideset{}{^*}\sum_{(n,m)\in L_k}(nm)^c a_nb_m\Big(\frac{n}{m}\Big) \\ +&\sideset{}{^*}\sum_{(n,m)\in L}(nm)^c a_nb_m\Big(\frac{n}{m}\Big).
\end{align*}
We conclude by bounding the second and third sums trivially. For the second we use the triangle inequality and then expand the sum to $I_k\times J'_k$:
\[
\sum_{k\in H}\;\sideset{}{^*}\sum_{(n,m)\in L_k}(nm)^c a_nb_m\Big(\frac{n}{m}\Big) \ll \sum_{k\in H}\;\sideset{}{^*}\sum_{n\in I_k}\;\sideset{}{^*}\sum_{m\in J'_k} T^c.
\]
Note that $\lvert I_k\times J'_k\rvert \leq \frac{T}{k(k+1)}$. Summing this over $k>z^{1/2}$ gives
\[
\sum_{k\in H}\;\sideset{}{^*}\sum_{(n,m)\in L_k}(nm)^c a_nb_m\Big(\frac{n}{m}\Big) \ll \frac{T^{1+c}}{z^{1/2}}.
\]
For the leftovers we use the triangle inequality again and expand the double sum to the region $((z,z+z^{1/2}]\times(z,\frac{T}{z}]\cup (T^{\delta}-z^{1/2},T^{\delta}]\times (z,\frac{T}{T^{\delta}-z^{1/2}}])\cap\mathbb{N}^2$. The number of integer pairs in this region is $\ll \frac{T}{z^{1/2}} + T^{1-\delta/2}\ll \frac{T}{z^{1/2}}$ using the assumption $z< T^{\delta}$. Thus we obtain the bound
\[
\sideset{}{^*}\sum_{(n,m)\in L}(nm)^c a_nb_m\Big(\frac{n}{m}\Big)\ll \frac{T^{1+c}}{z^{1/2}}.
\]
Overall, this gives the expression
\[
\sideset{}{^*}\sum_{z < n \leq T^{\delta}\;}\sideset{}{^*}\sum_{\;z < m\leq \frac{T}{n}} (nm)^c a_n b_m \Big(\frac{n}{m}\Big) = \sum_{k\in H}\;\;\sideset{}{^*}\sum_{n\in I_k}\;\;\sideset{}{^*}\sum_{m\in J_k} (nm)^c a_n b_m \Big(\frac{n}{m}\Big) + O\left(\frac{T^{1+c}}{z^{1/2}}\right).
\]
\end{proof}

We now complete our proof of Theorem \ref{general2}. This Theorem follows directly from Theorem \ref{general1} whenever $z\geq(\log T)^{24}$: in these cases, $\frac{(\log T)^2}{z^{1/12}}=O(1)$, so we obtain
\[
\sideset{}{^*}\sum_{\substack{z< n,m\leq T\\nm\leq T}}(nm)^c a_nb_m\Big(\frac{n}{m}\Big) \ll_c \frac{T^{1+c}(\log T)}{z^{1/4}}.
\]
We are left with the case where $z<(\log T)^{24}$, for which we aim to apply \eqref{Elliot} and Lemma \ref{cover}. However, in order for \eqref{Elliot} to be effective, we cannot allow $N=\max(I_k)$ to exceed $M^{1/2}=(\max(J_k))^{1/2}$ in any of our covering rectangles, as then the $N^2(\log N)$ term in \eqref{Elliot} would dominate the $M$ term, and may lead to bounds which are too large for our purposes. To avoid this we split the hyperbolic region as before:
\[
\sideset{}{^*}\sum_{\substack{z< n,m\leq T\\nm\leq T}}(nm)^c a_nb_m\Big(\frac{n}{m}\Big) = R_1(T) + R_2(T) + R_3(T) - R_4(T)
\]
where
\begin{align*}
    R_1(T)&=\sideset{}{^*}\sum_{\substack{T^{1/4} < n,m \leq T\\nm\leq T}} (nm)^c a_n b_m \Big(\frac{n}{m}\Big);\\
    R_2(T)&=\sideset{}{^*}\sum_{z < n \leq T^{1/4}\;}\sideset{}{^*}\sum_{\;z < m\leq \frac{T}{n}} (nm)^c a_n b_m \Big(\frac{n}{m}\Big);\\
    R_3(T)&= \sideset{}{^*}\sum_{z < m \leq T^{1/4}\;}\sideset{}{^*}\sum_{\;z < n\leq \frac{T}{m}} (nm)^c a_n b_m \Big(\frac{n}{m}\Big);\\
    R_4(T)&=\sideset{}{^*}\sum_{z < n \leq T^{1/4}\;}\sideset{}{^*}\sum_{\;z < m\leq T^{1/4}} (nm)^c a_n b_m \Big(\frac{n}{m}\Big).
\end{align*}
This splitting allows us to apply Lemma \ref{cover} to $R_2(T)$ (and $R_3(T)$) and obtain integer intervals $I_k$ whose maximums do not get too large, therefore allowing us to apply \eqref{Elliot} effectively. First, we bound $R_1(T)$ and $R_4(T)$. For $R_1(T)$ we use Theorem \ref{general1} with $z_1=T^{1/4}$:
\[
R_1(T) \ll_c \frac{T^{1+c}(\log T)^3}{z_1^{1/3}} = T^{11/12+c}(\log T)^3.
\]
Next, we use \eqref{HB} with $N=M=T^{1/4}$ to deal with $R_4(T)$. Multiplying by $1=\frac{T^c}{T^c}$ and setting $\widetilde{a}_n = \frac{n^c}{T^{c/2}}a_n$, $\widetilde{b}_m = \frac{m^c}{T^{c/2}}b_m$, we have $\lvert\widetilde{a}_n\rvert$, $\lvert\widetilde{b}_m\rvert\leq 1$. Then applying \eqref{HB} gives
\[
R_4(T) \ll_{\epsilon} T^{3/8+c+\epsilon},
\]
which is sufficient by choosing $\epsilon < 13/24$, as it may then be absorbed into the bound for $R_1(T)$.\\

We are left with $R_2(T)$ and $R_3(T)$. Note that these sums are symmetrically equivalent using the same argument as that of $S_2(T)$ and $S_3(T)$ in Section \ref{gen1proof}. Thus we only need to deal with $R_2(T)$. For this we use the covering Lemma \ref{cover} with $\delta = 1/4$:
\begin{align*}
    R_2(T) = \sum_{k\in H}\;\;\sideset{}{^*}\sum_{n\in I_k\;}\;\;\sideset{}{^*}\sum_{\;m\in J_k} (nm)^c a_n b_m \Big(\frac{n}{m}\Big) + O\left(\frac{T^{1+c}}{z^{1/2}}\right),
\end{align*}
where
\begin{align*}
H&=\left\{k\in\mathbb{N}:z^{1/2}\leq k\leq \frac{T^{1/4}}{z^{1/2}}-1\right\};\\
I_k &= \left\{n\in\mathbb{N}: z^{1/2}k<n\leq z^{1/2}(k+1)\right\};\\
J_k &= \left\{m\in\mathbb{N}: z<m\leq \frac{T}{z^{1/2}(k+1)}\right\}.
\end{align*}
To deal with this sum, we will consider the sum over $n$ and $m$ for a fixed $k$. First deal with the power term:
\begin{align*}
\sideset{}{^*}\sum_{n\in I_k}\;\;\sideset{}{^*}\sum_{m\in J_k} (nm)^c a_n b_m \Big(\frac{n}{m}\Big) &= T^c\sideset{}{^*}\sum_{m\in J_k}\;\;\sideset{}{^*}\sum_{n\in I_k} \frac{n^c}{ z^{c/2}(k+1)^c}\frac{m^cz^{c/2}(k+1)^c}{T^c}a_n b_m \Big(\frac{n}{m}\Big).
\end{align*}
Now $\frac{n}{z^{1/2}(k+1)}$, $\frac{mz^{1/2}(k+1)}{T}\leq 1$, so that we may define the sequences $\widetilde{a}_n = \frac{n^c}{z^{c/2}(k+1)^c}a_n$ in addition to $\widetilde{b}_m = \frac{m^cz^{c/2}(k+1)^c}{T^c}b_m$ which satisfy the condition $\lvert\widetilde{a}_n\rvert$, $\lvert\widetilde{b}_m\rvert\leq 1$. Finally we apply the Cauchy--Schwarz inequality and \eqref{Elliot} with $M= \max(J_k) \leq \frac{T}{z^{1/2}(k+1)}$ and $N = \max(I_k) \leq z^{1/2}(k+1)$:
\begin{align*}
\sideset{}{^*}\sum_{n\in I_k}\;\;\sideset{}{^*}\sum_{m\in J_k} (nm)^c a_n b_m \Big(\frac{n}{m}\Big) &= T^c\sideset{}{^*}\sum_{m\in J_k}\;\;\sideset{}{^*}\sum_{n\in I_k} \widetilde{a}_n \widetilde{b}_m \Big(\frac{n}{m}\Big)\\
\\ &\ll T^c\lvert J_k\rvert^{1/2} \Big(\sideset{}{^*}\sum_{m\in J_k} \Big\lvert\sideset{}{^*}\sum_{n\in I_k} \widetilde{a}_n \Big(\frac{n}{m}\Big)\Big\rvert^2\Big)^{1/2}\\
&\ll T^c \Big(\frac{T}{z^{1/2}(k+1)}\Big)^{1/2}\Big(\frac{T}{z^{1/2}(k+1)}\Big)^{1/2}\lvert I_k\rvert^{1/2}\\
&\ll \frac{T^{1+c}}{z^{1/4}(k+1)}
\end{align*}
where we used the fact that $k^2 z (\log(kz^{1/2})) \ll T^{1/2}(\log T)$, while $\frac{T}{z^{1/2}(k+1)}\geq T^{3/4}$ to simplify the application of \eqref{Elliot}. Summing this bound over the given $k$ introduces a logarithmic term, and so
\begin{align*}
    R_2(T), R_3(T) &\ll \frac{T^{1+c}(\log T)}{z^{1/4}}.\\
\end{align*}
Combining all the bounds we get:
\[
\sideset{}{^*}\sum_{\substack{z< n,m\leq T\\nm\leq T}}(nm)^c a_nb_m\Big(\frac{n}{m}\Big) \ll_c \frac{T^{1+c}(\log T)}{z^{1/4}} + T^{11/12+c}(\log T)^3\ll_c \frac{T^{1+c}(\log T)}{z^{1/4}}
\]
(since $z<(\log T)^{24}$) as required.

\section{Proof of Theorem \ref{asymptotic}}
To begin we once more cut the hyperbolic region into regions depending on the sizes of each variable. We write
\[
\sum_{\substack{1\leq n,m\leq T\\2\nmid nm\\nm\leq T}} \Big(\frac{n}{m}\Big) = N_1(T) + N_2(T) - N_3(T),
\]
where
\begin{align*}
    N_1(T) &= \sum_{\substack{1\leq n\leq T^{1/2}\\2\nmid n}}\sum_{\substack{1\leq m\leq T/n\\2\nmid m}} \Big(\frac{n}{m}\Big);\\
    N_2(T) &= \sum_{\substack{1\leq m\leq T^{1/2}\\2\nmid m}}\sum_{\substack{1\leq n\leq T/m\\2\nmid n}} \Big(\frac{n}{m}\Big);\\
    N_3(T) &= \sum_{\substack{1\leq n\leq T^{1/2}\\2\nmid n}}\sum_{\substack{1\leq m\leq T^{1/2}\\2\nmid m}} \Big(\frac{n}{m}\Big).
\end{align*}
Let us first deal with $N_1(T)$. We begin by separating the square values of $n$:
\[
N_1(T) = \sum_{\substack{1\leq n\leq T^{1/2}\\2\nmid n\\n = \square}}\sum_{\substack{1\leq m\leq T/n\\2\nmid m}} \Big(\frac{n}{m}\Big) + \sum_{\substack{1\leq n\leq T^{1/2}\\2\nmid n\\n\neq\square}}\sum_{\substack{1\leq m\leq T/n\\2\nmid m}} \Big(\frac{n}{m}\Big).
\]
To deal with the second of these sums we use the P\'olya--Vinogradov inequality for the sum over $m$, and then sum over $1\leq n\leq T^{1/2}$. Thus the second sum is $O(T^{3/4}(\log T))$. For the first sum, we note that since $n$ is a square the Jacobi symbol is the trivial character modulo $n$. Thus
\[
\Big(\frac{n}{m}\Big) =
\begin{cases}
1\; \text{if} \; \gcd(n,m)=1,\\
0\; \text{if} \; \gcd(n,m)>1.
\end{cases}
\]
It is well-known that for a fixed odd $n$, the number of odd $1\leq m\leq \frac{T}{n}$ co-prime to $n$ is given by
\[
\frac{T}{2n}\cdot\frac{\phi(n)}{n} + O(n^{\epsilon})
\]
for any $\epsilon>0$. Summing this error over the square values of $n$ less than $T^{1/2}$ we will obtain an error of size $O(T^{1/4+\epsilon})$, which is satisfactory. For the main term we use the change of variables $n=k^2$:
\begin{align*}
    \sum_{\substack{1\leq n\leq T^{1/2}\\2\nmid n\\n=\square}} \frac{T}{2n}\cdot\frac{\phi(n)}{n} &=\frac{T}{2} \sum_{\substack{1\leq k\leq T^{1/4}\\2\nmid k}} \frac{\phi(k^2)}{k^4} =  \left(\sum_{\substack{k=1\\2\nmid k}}^{\infty} \frac{\phi(k^2)}{k^4}\right)\frac{T}{2} + O(T^{3/4}).
\end{align*}
Noting that $N_2(T)$ may be dealt with using the same methods we obtain
\[
N_1(T)+N_2(T) = \left(\sum_{\substack{k=1\\2\nmid k}}^{\infty} \frac{\phi(k^2)}{k^4}\right)T + O(T^{3/4}(\log T)).
\]
Using Theorem $1$ of \cite{CFS} with $X=Y=T^{1/2}$ we obtain, $N_3(T) \ll T^{3/4}$. Thus we have
\[
\sum_{\substack{1\leq n,m\leq T\\2\nmid nm\\nm\leq T}} \Big(\frac{n}{m}\Big) = \left(\sum_{\substack{k=1\\2\nmid k}}^{\infty} \frac{\phi(k^2)}{k^4}\right)T + O(T^{3/4}(\log T)).
\]
Lastly we evaluate the constant. To do this, let $g(n) =\mathds{1}_{\textrm{odd}}(n)\mathds{1}_{\square}(n)\phi(n)$ where $\mathds{1}_{\textrm{odd}}$ and $\mathds{1}_{\square}$ are the indicator functions for odd numbers and squares respectively. We will consider the Dirichlet series and Euler product of this multiplicative function:
\begin{align*}
\sum_{\substack{k=1}}^{\infty} \frac{g(n)}{n^s} &= \prod_{p}\left(1+\sum_{m=1}^{\infty}\frac{g(p)}{p^{ms}}\right)\\
&= \prod_{p\neq 2}\left(1+\sum_{m=1}^{\infty}\frac{\phi(p^{2m})}{p^{2ms}}\right)\\
&= \prod_{p\neq 2}\left(1+\frac{(p-1)}{p}\sum_{m=1}^{\infty}\frac{1}{p^{2m(s-1)}}\right)\\
&= \prod_{p\neq 2}\left(\frac{1-1/p^{(2s-1)}}{1-1/p^{(2s-2)}}\right)\\
&= \frac{1-2^{(2-2s)}}{1-2^{(1-2s)}}\frac{\zeta(2s-2)}{\zeta(2s-1)},
\end{align*}
where $\zeta$ is the Riemann-zeta function. By taking $s=2$ we obtain the equality
\[
\sum_{\substack{k=1\\2\nmid k}}^{\infty} \frac{\phi(k^2)}{k^4} = \frac{6\zeta(2)}{7\zeta(3)}
\]
as required.

\end{document}